\documentclass[12pt,a4paper]{article}

\usepackage[T1]{fontenc}
\usepackage[utf8]{inputenc}
\usepackage{authblk}
\usepackage{amsfonts}
\usepackage{mathtools}
\usepackage{slashed}
\usepackage{amsmath,amssymb}
\usepackage{graphicx,xcolor}
\usepackage{cite}
\usepackage{hyperref}

\definecolor{refcol}{rgb}{0.9,0.1,0.1}
\hypersetup{colorlinks=true,linkcolor=blue,citecolor=refcol,urlcolor=cyan,linktocpage}
  
\usepackage[capitalize]{cleveref}
\usepackage{tikz}
\usetikzlibrary{shapes}
\usetikzlibrary{plotmarks}

\usepackage{setspace}

\usepackage[left=2.5cm,right=2.5cm,top=2.5cm,bottom=2.5cm]{geometry}
\allowdisplaybreaks
\newcommand{\be}{\begin{equation}}
\newcommand{\ee}{\end{equation}}
\newcommand{\bea}{\begin{eqnarray}}
\newcommand{\eea}{\end{eqnarray}}

\def\XXint#1#2#3{{\setbox0=\hbox{$#1{#2#3}{\int}$ }
\vcenter{\hbox{$#2#3$ }}\kern-.6\wd0}}

\newcommand{\Qp}{\mathbb{Q}_p}
\newcommand{\Zp}{\mathbb{Z}_p}

\begin{document}
\begin{titlepage}

\title{
{\huge\bf Pseudodifferential Operators on $\Qp$}\\ 
{\Huge\bf and $L$-Series}
}

\bigskip\bigskip\bigskip\bigskip\bigskip

\author{{\bf Parikshit Dutta}${}^{1}$\thanks{{\tt parikshitdutta@yahoo.co.in}}\hspace{8pt} {\bf and} \hspace{2pt}  
                    {\bf Debashis Ghoshal}${}^2$\thanks{{\tt dghoshal@mail.jnu.ac.in}} \\  
\hfill\\
${}^1${\it Asutosh College, 92 Shyama Prasad Mukherjee Road,}\\
{\it Kolkata 700026, India}
\hfill\\              
${}^2${\it School of Physical Sciences, Jawaharlal Nehru University,}\\
{\it New Delhi 110067, India}
}

\date{%
%
\bigskip\bigskip
\begin{quote}
\centerline{{\bf Abstract}}
{\small
We define a family of pseudodifferential operators on the Hilbert space $L^2(\Qp)$ of complex valued square-integrable 
functions on the $p$-adic number field $\Qp$. The Riemann zeta-function and the related Dirichlet $L$-functions can be 
expressed as a trace of these operators on a subspace of $L^2(\Qp)$. We also extend this to the $L$-functions associated 
with modular (cusp) forms. Wavelets on $\Qp$ are common sets of eigenfunctions of these operators. 
}
\end{quote}
}

\bigskip

\end{titlepage}
\thispagestyle{empty}\maketitle\vfill \eject

\tableofcontents

\section{Introduction}\label{sec:Intro} 
The Riemann zeta function\cite{hmedwards} has a product representation in terms of prime numbers that allows one to connect 
to the $p$-adic number field $\Qp$. This prompted us (in a larger collaboration) to construct a unitary matrix model (UMM) for 
the former by combining UMMs for the local factors corresponding to each prime\cite{Chattopadhyay:2018bzs}. The parameters 
that define the UMM for the Riemann zeta function is divergent, however, after a renormalisation, resultant UMM agrees with that
constructed directly in Ref.\cite{Dutta2016byx}. Moreover, we show that the partition function of the UMM can be written as a trace 
of the (generalised) Vladimirov derivative \cite{VVZ1994p} acting on a (subspace) of complex valued square-integrable functions 
on the $p$-adic number field $\Qp$. The partition function is essentially  the Riemann zeta function in this approach. 

While integration of complex valued valued functions on $\Qp$ is quite straightforward, the totally disconnected topology of 
$\Qp$ does not allow for a na\"{i}ve notion of a derivative. The Valdimirov derivative is, therefore, defined through an integral 
kernel. Moreover, Kozyrev\cite{Kozyrev:2001} (see also \cite{Albeverio:2006,Albeverio:2009,Albeverio:2011}) constructed 
complex valued wavelets on $\Qp$ and demonstrated that these are eigenfunctions of the generalised Vladimirov derivative.
These wavelets are rather like the complex valued generalised Haar wavelets on the real line $\mathbb{R}$. The eigenvalues
are related to the scaling property of the wavelet. Enhancement of the manifest affine symmetry of the wavelets were studied 
in \cite{Dutta:2018qrv}.     

The zeta function of Riemann is but one member of an infinite family of Dirichlet $L$-functions, which are defined using 
arithmetic functions called Dirichlet characters. These are multiplicative characters, hence, all Dirichlet $L$-functions admit
Euler product representations over prime numbers. It is, therefore, imperative that we generalise the notion of Vladimirov 
derivative further and seek pseudodifferential operators, the eigenvalues of which will involve the Dirichlet character. 

Indeed, it turns out that one can do more. There are $L$-functions associated with modular forms, which we refer to as modular 
$L$-function for brevity. These functions too have both Dirichlet series and Euler product representations. We are able to define 
a pair of pseudodifferential operators corresponding to each prime factor. We discuss the eigenvalues and eigenfunctions (which 
are the Kozyrev wavelets again) of these generalised pseudodifferential operators. The $L$-functions can be written in terms of the 
trace of these operators on a subspace of the Hilbert space spanned by the wavelets. We also comment on the realisation of the
Hecke operarors as traces of these operators after conjugation by raising-lowering operators.
 

\section{Dirichlet characters and Dirichlet $L$-functions}\label{sec:DChar}
For a positive integer $k$, the \emph{Dirichlet character modulo $k$} is a map $\chi_k : \mathbb{Z} \rightarrow
\mathbb{C}$ such that
\begin{enumerate}
\item
For all $m_1, m_2 \in \mathbb{Z}$, $\chi_k(m_1 m_2) = \chi_k(m_1) \chi_k(m_2)$
\item
$\chi_k(m_1) = \chi_k(m_2)$ if $m_1 \equiv m_2$ (mod $k$)
\item
$\chi_k(m) \ne 0$ if and only if $m$ is relatively prime to $k$
\end{enumerate}
There is a \emph{trivial character} that assigns the value 1 to all integers (including 0). This corresponds to $k=1$. 

More precisely, the multiplicative group $G(k)=\left(\mathbb{Z}/k\mathbb{Z}\right)^*$ consisting of the \emph{invertible} elements 
of $\mathbb{Z}/k\mathbb{Z}$ is an abelian group. An element $\chi_k\in \mathrm{Hom}(G(k),\mathbb{C}^*)$ is called a 
character modulo $k$. It is a $\mathbb{C}^*$-valued function on the set of integers relatively prime to $k$, such that the
property (i) above holds. It is convenient and often conventional to \emph{extend} a character to all $\mathbb{Z}$ by setting
$\chi_k(m)=0$ for all $m$ which are not coprime to $k$ \cite{Serre:Course}.    

The periodicity $k$ of the Dirichlet character does not specify the function completely. For any $k$ there are $\varphi(k)$
number\footnote{For any positive integer $n$, the Euler totient function $\varphi(n)$ counts the number of positive  integers 
$1< m <n$ that are relatively prime to $n$.} of inequivalent characters. There is a \emph{principal character} $\chi_{k,0}$  that 
assumes the value 1 for arguments coprime to $k$ and vanishes otherwise. If all the values of a character is real, it is called a 
\emph{real character}. Clearly, if $k$ is a prime number, there are $(k -1)$ inequivalent characters, the values of which are 
$\varphi(k)$-th roots of unity. Of these the character $\chi_{k,0}(n)$ vanishes for all integers $n\equiv 0$ (mod $k$) and 1 
otherwise. If $k$ is not a prime, the character vanishes for all integers that share common (prime) factors with it.    

The Dirichlet $L$-series corresponding to a Dirichlet character $\chi$, is the infinite series 
\begin{equation}
L(s,\chi) = \sum_{n=1}^\infty \frac{\chi(n)}{n^s}\label{LSeries}
\end{equation}
which is convergent for $\mathrm{Re }(s) > 1$. Its analytic continuation in the complex $s$-plane defines a meromorphic 
function called the Dirichlet $L$-function, also denoted by the same symbol. If $\chi$ is chosen to be the trivial character, 
one gets the Riemann zeta function by the analytic continuation of
\begin{equation}
\zeta(s) = \sum_{n=1}^\infty \frac{1}{n^s} =  \prod_{p\,\in\,{\mathrm{primes}}} \frac{1}{\left(1 - p^{-s}\right)} 
\label{EulerSum}
\end{equation}
where, the second expression is the Euler product form. Unlike the Riemann zeta function, the Dirichlet functions $L(s,\chi)$ 
defined by the analytic continuation of \cref{LSeries} is an entire function, except when $\chi$ is the principal character (or the
trivial one), in which case there is a simple pole at $s=1$. Moreover, they have a set of trivial zeroes at negative even
or odd integers (depending on whether $\chi(-1)=\pm1$) and a set of non-trivial zeroes which must lie symmetrically about the 
\emph{critical line} $\mathrm{Re }(s) = \frac{1}{2}$ in the critical strip $0 < \mathrm{Re }(s) < 1$. The assertion that the non-trivial
zeroes are exactly on the critical line is the \emph{generalised Riemann hyposthesis}.

Since the Dirichlet characters are multiplicative, the $L$-series \cref{LSeries} can be written as an infinite product over
prime numbers
\begin{equation}
L(s,\chi)  = \prod_{p\,\in\,{\mathrm{primes}}} \frac{1}{\left(1 - \chi(p) p^{-s}\right)}\label{EulerProd}
\end{equation}
in the region of its convergence. We shall be interested in studying the \emph{local} factor $\left(1 - \chi(p) p^{-s}\right)^{-1}$
at the prime $p$, which is the sum of an infinite geometric series.

\bigskip

Thanks to the multiplicative nature of a Dirichlet character, $\chi_k(m^n) = \left(\chi_k(m)\right)^n$ for $n\ge 0$. In other words,
Dirichlet characters of a positive integer power of an integer is well defined. However, since negative powers of an integer is not 
an integer, $\chi_k(m^{-n})$ is not defined. In order to define the pseudodifferential operators we are interested in, however, we 
shall need to \emph{extend} the notion of a Dirichlet character to negative integer powers as well. Therefore, let us \emph{define}
the function $\mathfrak{x}_k(m^{n})$ as 
\begin{equation}
\mathfrak{x}_k(m^{n}) = \left\{\begin{array}{ll}
 \left(\chi_k(m)\right)^{n} & \text{ if } \: \chi_k(m) \ne 0\:\text { and }\, n\ge 0\\
 \left(\chi_k(m)\right)^{n} = \chi^*_k(m^{-n})& \text{ if } \: \chi_k(m) \ne 0\:\text { and }\, n < 0\\
 0 & \text{ if } \: \chi_k(m) = 0\:\text { for all }\, n\ne 0 \end{array} 
  \right. \label{XDirChar}
\end{equation}
Thus, if the Dirichlet character of an integer is zero, we take all its integer powers (positive or negative) to be zero as well. In this paper 
we  shall refer to $\mathfrak{x}_k$ as an \emph{extended Dirichlet character}. In fact, for our purposes, it will be sufficient to define the 
extended character for integer powers of prime numbers\footnote{The extension we need is really quite minimal. The original Dirichlet 
character is a map $\chi: \mathbb{Z} \rightarrow \mathbb{C}$. This of course makes sense as a map $\chi\big|_{(p)} : p^{\mathbb{N}} 
\rightarrow \mathbb{C}$. We need to extend this restricted form to $\mathfrak{x} : p^{\mathbb{Z}} \rightarrow \mathbb{C}$. Moreover, 
only the non-zero values of the characters are precisely defined, the rest was by extension. Therefore, our definition should also be 
thought of as an extension in the same spirit.} only. That is because, we need this extension not on all $p$-adic numbers, but only on 
their norms\footnote{It should also be noted that the extended character we need is really a combination of the `norm function' on 
$\Qp$ ($| \cdot|_p : \Qp \rightarrow \mathbb{R}$) and an extended character on rational numbers. It may be possible to define an 
extended character on $\Qp$ itself, however, the present construction serves our purpose.}. Consequently, if $k$ contains $p$ (of 
$\Qp$) as one of its factors, the extended character vanishes for all $\xi\in\Qp$. On the other hand, if $k$ and $p$ are relatively prime, 
its contribution is a constant phase on a subset of $\Qp$ with a fixed norm. 

\section{Twisted $p$-adic Gelfand-Graev-Tate Gamma functions}\label{sec:GGTGamma}
With the help of the extended Dirichlet character $\mathfrak{x}$ in \cref{XDirChar}, let us define the Gelfand-Graev-Tate gamma 
function $\Gamma_\mathfrak{x}$, twisted by $\mathfrak{x}_k$ as
\begin{equation}
\Gamma_\mathfrak{x} (s) = \int_{\Qp^\times} \frac{d\xi}{|\xi|_p} \, e^{2\pi i\xi} |\xi|_p^s\, \mathfrak{x}_k\!\left(|\xi|_p^{-1}\right)  
= \int_{\Qp} d\xi \, e^{2\pi i\xi} |\xi|_p^{s-1}\, \mathfrak{x}_k\!\left(|\xi|_p^{-1}\right) , \quad p\nmid k
\label{TwistedGGTGamma}
\end{equation}
If $p$ divides $k$, then the function vanishes identically. The definitions above as well as in \cref{GammA} to appear below are 
consistent with the definition of a generalised gamma function with a multiplicative character \cite{Gelfand1968representation}
(similar functions have been used in \cite{Ghoshal:2004ay,Gubser:2017qed}). We shall now evaluate this integral. 

Recall the `Laurent expansion' of a $p$-adic number $\xi = p^n\left(\xi_0 + \xi_1 p + \xi_2 p^2 + \cdots\right)$, $p\in 
\mathbb{Z}$, $\xi_m\in\{0,1,\dots,p-1\}$ but $\xi_0\ne 0$. As usual, we split the integral into three parts corresponding to 
the regions 
\begin{enumerate}
\item[(i)] 
$|\xi|_p < 1$. Hence $e^{2\pi i \xi} = 1$. The region can be further divided into `circles' $C_n = \left\{\xi\; :\; |\xi|_p = p^{-n}\right\}$, 
$n=1,2,\cdots$, the measure of which is $\frac{p-1}{p} p^{-n}$. Thus the contribution to the integral from this region is
\[
\frac{p-1}{p} \sum_{n=1}^\infty p^{-n}\, p^{-n(s-1)}\, \mathfrak{x}_k(p^{n}) 
= \frac{p-1}{p} \sum_{n=1}^\infty \left( \chi_k(p) p^{-s}\right)^{n}  
= \frac{p-1}{p} \frac{\chi_k(p)}{p^s - \chi_k(p)}  
\]
for non-zero $\chi_k(p)$.
\item[(ii)] $|\xi|_p = 1$. Hence $\mathfrak{x}_k(|\xi_p|)=\chi_k(|\xi_p|)=1$ and $e^{2\pi i\xi}=1$. The measure of this region 
is $\frac{p-1}{p}$. So the contribution to the integral is
$\frac{p-1}{p}$. 
\item[(iii)] $|\xi|_p > 1$. It consists of `circles'  $C_n = \left\{\xi\; :\; |\xi|_p = p^{n}\right\}$, $n=1,2,\cdots$. In the circle $C_n$,
$\mathfrak{x}_k(|\xi|_p) = \mathfrak{x}_k(p^{-n}) = \left(\chi_k(p)\right)^{-n}$ and the measure of $C_n$ is $\frac{p-1}{p} p^{n}$. The 
circle $C_n$ is divided into $(p-1)$ subsets $C_{n;1}, C_{n;2}, \cdots, C_{n;p-1}$ of equal measures. This, in turn is further divided 
into subsets such that
\[
C_{n} = \displaystyle{\bigcup_{m_0=1}^{p-1}} C_{n;m_0} = \displaystyle{\bigcup_{m_0,m_1=0}^{p-1}} C_{n;m_0,m_1} = 
\cdots = \displaystyle{\bigcup_{m_0,\cdots,m_{n-1}=0}^{p-1}} \!\!\!\!C_{n;m_0,m_1,\cdots,m_{n-1}}
\]
In $C_{n;m_0,m_1,\cdots,m_{n-1}}$, the exponential $e^{2\pi i\xi} = \omega_{p^n}^{\xi_0} \omega_{p^{n-1}}^{\xi_1} \cdots 
\omega_p^{\xi_{n-1}}$. What is important is the last factor, which is a $p$-th roots of unity. When integrated, i.e., summed 
over all values $m_{n-1}=0,1,\cdots,p-1$, the sum over roots of units vanish (all subsets have equal measure). Therefore, the 
contribution to the integral from $C_n$ is zero, except for $n=1$. In this case, the root 1 is missing, so the contribution to the 
integral is
$\frac{p-1}{p} \frac{p}{p-1} \left( \omega_p + \omega_p^2 + \cdots + \omega_p^{p-1}\right) p^{s-1} \mathfrak{x}_k(p^{-1})
= - \displaystyle{\frac{p^{s-1}}{\chi_k(p)}}$,  
if $\chi(p)\ne 0$, and 0 otherwise. 
\end{enumerate}
Combining the three contributions
\begin{equation}
\Gamma_\mathfrak{x} (s) = \frac{p^s\, \left(\chi_k(p) - p^{s-1}\right)}{\chi_k(p) \left(p^s - \chi_k(p)\right)}
= \frac{\chi_k(p) -  p^{s-1}}{\chi_k(p)\left(1 - \chi_k(p) p^{-s}\right) } 
\label{GGTvalue}
\end{equation}
For the trivial character this reduces to $\Gamma (s) = \displaystyle\frac{1 - p^{s-1}}{1 - p^{-s}}$, the standard Gelfand-Graev-Tate 
gamma function. In the general case, \cref{TwistedGGTGamma} further generalises the gamma function defined with multiplicative 
characters in, e.g, Ref.\cite{Ghoshal:2004ay,Gubser:2017qed}. 

\section{Vladimirov derivative twisted by a character}\label{sec:ChiVlad}
Let us recall that the generalised Vladimirov derivative, defined as an integral kernel
\begin{equation}\label{VladD}
D^\alpha g(\xi) = \frac{1}{\Gamma(-\alpha)}\, \int d\xi'\, \frac{g(\xi') - g(\xi)}{|\xi'-\xi|_p^{\alpha+1}}
\end{equation}
is defined for any $\alpha\in\mathbb{C}$ by analytic continuation. The eigenfunctions of this pseudo-differential operator

\begin{equation}
D^\alpha \psi_{n,m,j}^{(p)} (\xi) = p^{\alpha(1-n)} \psi_{n,m,j}^{(p)} (\xi) \label{VDonKozy}
\end{equation}
with eigenvalue $p^{\alpha(1-n)}$ are the Kozyrev wavelets
\begin{equation}
\psi_{n,m,j}^{(p)} (\xi) = p^{-\frac{n}{2}} e^{2\pi i j p^{n-1} \xi}\, \Omega^{(p)}\! \left( p^n \xi - m \right), 
\quad \xi \in \Qp 
\label{pWavelet} 
\end{equation}
for $n\in\mathbb{Z}$, $m \in \Qp/\Zp$ and $j \in \{1,2,3,\cdots,p-1\}$. These are Bruhat-Schwartz functions on $\Qp$, consisting of 
the indicator function\footnote{The indicator function $\Omega^{(p)}(\xi) = 1$ for $|\xi|_p \le 1$, and 0 otherwise.} $\Omega^{(p)}(\xi)$
and an additive character, the exponential function. 

As expected the eigenvalues depend only on the quantum number related to scaling $n$ (and not on those related to translation 
and phase). Since we are primarily interested in the eigenvalues, we may as well restrict our attention to the set of eigenfunctions 
$\psi_{n,0,1}^{(p)} (\xi) = p^{-\frac{n}{2}} e^{2\pi i p^{n-1} \xi}\, \Omega^{(p)}\! \left( p^n \xi \right)$, which we shall often denote 
simply as $\left| 1 - n\right\rangle$ labelled by the eigenvalues. It is natural to define the raising-lowering operators $a_\pm$ 
\begin{equation}
a_\pm\psi_{n,0,1}(\xi) = \psi_{n\pm 1,0,1}(\xi), \qquad  a_\pm | n \rangle = | n \mp 1\rangle
\label{raise:lower}
\end{equation}
that changes the scaling quantum number by one. These operators, together with $\log_p D = 
\displaystyle{\lim_{\alpha\to 0}}\frac{D^\alpha - 1}{\alpha\ln p}$, (formally) generate a large symmetry of the wavelets 
\cite{Dutta:2018qrv}. 

For now it suffices to note that
\begin{equation}
D^\alpha \left| n\right\rangle = p^{\alpha n}  \left| n\right\rangle \label{VladDonKet}
\end{equation}
which allows us to write the local factor $\left(1-p^{-s}\right)^{-1}$ in the Euler product form of the Riemann zeta function 
\cref{EulerSum} as 
\begin{equation}
\zeta_p(s) = \frac{1}{\left(1-p^{-s}\right)} = \sum_{m=0}^\infty p^{-sm} = \sum_{m=0}^\infty \left\langle m \right| D^{-s} 
\left| m\right\rangle = \mathrm{Tr}_{\mathcal{H}_-} D^{-s} \label{locZTrace}
\end{equation}
i.e., a trace of the generalised Vladimirov derivative $D^{-s}$ over a subspace $\mathcal{H}_{-}$ spanned by the wavelets 
$\left\{|m\rangle \sim \psi_{1-m,0,1}\,:\, m=1,2,\cdots\right\}$. In other words, the subspace $\mathcal{H}_{-}$ is spanned by  
wavelets on $\Zp$ with compact support. 

\bigskip

Let us now define the generalised Vladimirov derivative, \emph{twisted by the character} $\mathfrak{x}$, as  
\begin{equation}\label{ChiVladD}
D_\mathfrak{x}^\alpha g(\xi) = \frac{1}{\Gamma_\mathfrak{x}(-\alpha)}\, \int_{\Qp} d\xi'\, \frac{g(\xi') - g(\xi)}{|\xi'-\xi|_p^{\alpha+1}}\, 
\mathfrak{x}_k(|\xi'-\xi|_p^{-1}) 
\end{equation}
Once again, this operator is meaningful if $p\nmid k$. If, on the other hand, $p$ divides $k$, the twisted character vanishes identically.
consequently the numerator on the RHS is zero. However, the twisted gamma function $\Gamma_{\mathfrak{x}}$ in the denominator 
also vanishes. In this case, we define $D_{\mathfrak{x}}=1$ if $(p,k)\ne 1$. 


We shall now show that the Kozyrev wavelets \cref{pWavelet} are also eigenfunctions of these operators. Actually due to the 
translation, scaling and the phase rotation properties, it is sufficient to establish that $f(\xi) = \psi_{0,0,1} = e^{2\pi i \xi/p} 
\Omega^{(p)} \left( \xi \right)$ is an eigenfunction. The rest follows trivially. The proof is in fact very similar to the case of 
the generalised Vladimirov derivative with trivial character, so we shall be brief. We consider two cases:
\begin{itemize}
\item 
Case (i) The set $|\xi|_p > 1$ is not in the support of $\Omega^{(p)} \left( \xi \right)$, therefore, $f(\xi) = 0$. Hence, the contribution
to the integral is only from the set $|\xi'|_p < 0$. By the triangle inequality 
$\mathfrak{x}_k(|\xi' - \xi|_p^{-1}) = \mathfrak{x}_k(|\xi|_p^{-1})$ and 
$|\xi' - \xi|_p^{\alpha+1} = |\xi|_p^{\alpha+1}$ does not have any dependence on the variable of integration. Thus the integral, which is 
reduced to $\displaystyle{\int_{|\xi'|_p\le 1} e^{2\pi i\xi'/p} d\xi'}$, is zero. 
\item 
Case (ii) Only those $\xi$ for which $|\xi|_p \le 1$ need to be considered.
The function $f(\xi)$ is supported on this set, hence, $\Omega^{(p)} \left( |\xi|_{p} \right) = 1$, therefore, the numerator of the integrand
involves $e^{2\pi i\xi'/p}\, \Omega^{(p)}\! \left( \xi'\right) - e^{2\pi i\xi/p} = e^{2\pi i\xi/p}\left(e^{2\pi i(\xi' - \xi)/p}\, \Omega^{(p)}\!\left( 
\xi'\right) - 1\right)$. The prefactor does not depend on the variable of integration. Due to the fact that every point in a $p$-adic 
disc could be considered to be its centre, we may write $\Omega^{(p)}\! \left( \xi'\right) = \Omega^{(p)}\! \left( \xi'-\xi\right)$ and change
the variable of integration to $z=\xi' - \xi$ by a translation. Now we can split the contribution to the integral from the three sets
$|z|_p < 1$, $|z|_p = 1$ and $|z|_p > 1$, of which the first one vanishes (since $e^{2\pi iz/p} = 1$ hence the integrand is zero).  
{}From the other two sets, we get
\begin{eqnarray}
D_\mathfrak{x}^\alpha f(\xi) 
&{=}& \frac{e^{2\pi i\xi/p}}{\Gamma_\mathfrak{x}(-\alpha)} \, \left[  \int_{|z|_p=1} \left(e^{2\pi iz/p}-1\right) dz - \int_{|z|_p > 1} \frac{\mathfrak{x}_k(|z|_p^{-1})}{|z|_p^{\alpha+1}}\, 
dz  \right] \nonumber\\
&{=}& \frac{e^{2\pi i\xi/p}}{\Gamma_\mathfrak{x}(-\alpha)} \, \left[\frac{1}{p} \left(\omega_p + \cdots + \omega_p^{p-1}\right) - \frac{(p-1)}{p} -\sum_{n=1}^\infty 
\frac{\mathfrak{x}_k(p^{-n})}{p^{n(\alpha+1)}}\, \frac{p^n(p-1)}{p} \right] \nonumber\\
&{=}& - \, \frac{e^{2\pi i\xi/p}\, \Omega^{(p)}\!\left(\xi \right)}{\Gamma_\mathfrak{x}(-\alpha)} \, \left( 1 + 
\frac{p-1}{p}\, \frac{1}{p^\alpha \chi_k(p) -1} \right) \nonumber\\
&{=}& p^\alpha \chi_k(p)\, f(\xi) \label{TVladEFn}
\end{eqnarray}
where we have used the fact that $\Omega^{(p)}\!\left(\xi \right) =1$.    
\end{itemize}
This proves that $f(\xi) = \psi_{0,0,1}(\xi) = e^{2\pi i \xi/p}\,\Omega^{(p)}\!\left( \xi \right)$ is an eigenfunction corresponding
to the eigenvalue $\chi_k(p)p^{\alpha}$. Hence,
\begin{align}
D^\alpha_\mathfrak{x} \psi_{n,m,j}^{(p)} (\xi) &= \left(\chi_k(p) p^{\alpha}\right)^{(1-n)} \psi_{n,m,j}^{(p)} (\xi) \nonumber\\
\text{or,}\qquad
D^\alpha_\mathfrak{x} \left| n \right\rangle &= \left(\chi_k(p) p^{\alpha}\right)^{n} \left| n \right\rangle
\label{VTwDonKozy}
\end{align}
where the latter expression is in the notation of \cref{VladDonKet}. We see that the family of pseudodifferential operators 
$D_\mathfrak{x}$, twisted by the multiplicative Dirichlet character $\mathfrak{x}$, all commute with each other, and their common 
eigenfunctions, with eigenvalues $\left(\chi_k(p) p^{\alpha}\right)^{(1-n)}$, are the Kozyrev wavelets on $\Qp$.

Repeating the arguments in \cref{locZTrace} \emph{mutatis mutandis}, we can write each local factor in the Dirichlet $L$-function 
as 
\begin{equation}
L_p(s,\chi) = \frac{1}{\left(1-\chi(p) p^{-s}\right)} = \sum_{m=0}^\infty \left(\chi(p) p^{-s}\right)^m = 
\sum_{m=0}^\infty \left\langle m \right| D_\mathfrak{x}^{-s} \left| m\right\rangle = \mathrm{Tr}_{\mathcal{H}_{-}} 
\left(D_\mathfrak{x}^{-s}\right) \label{locLTrace}
\end{equation}
i.e., a trace over a subspace $\mathcal{H}_{-}=L^2(\Zp)$ of $L^2(\Qp)$ spanned by the Kozyrev wavelets.

\section{$L$-functions of modular forms and pseudodifferential operators}\label{sec:VladDLFn}
A more general class of $L$-functions are those associated with a modular form $f$ 
\cite{Serre:Course,Koblitz:ECMF,Warner:Crash,MIT:MFLF}. The group $\mathrm{SL}(2,\mathbb{R})$ has a natural action on 
the upper half plane $\mathbb{H} = \{z \,:\, \mathrm{Im }(z) > 0\}$. Consider its discrete subgroup\footnote{More precisely, 
the relevant groups are the projective special linear groups PSL(2,$\mathbb{R}$) = SL(2,$\mathbb{R}/\{\pm\}$) and 
PSL(2,$\mathbb{Z}$) = SL(2,$\mathbb{Z})/\{\pm\}$.} SL(2,$\mathbb{Z})$ and its following \emph{congruence subgroups}
of finite indices: 
\begin{align}
\Gamma_0(N) &= \left\{\gamma\in \mathrm{SL}(2,\mathbb{Z})\, |\, c \equiv 0\: \text{mod}\: N\right\}\nonumber\\  
\Gamma_1(N) &= \left\{\gamma\in \mathrm{SL}(2,\mathbb{Z})\, |\, a,d \equiv 1\:\text{and}\: c \equiv 0\: \text{mod}\: N\right\}
\label{CongruentSubGp}\\
\Gamma(N) &= \left\{\gamma\in \mathrm{SL}(2,\mathbb{Z})\, |\, a,d \equiv 1\:\text{and}\: b,c \equiv 0\: \text{mod}\: N\right\}\nonumber
\end{align}     
The conditions are empty for $N=1$, hence $\Gamma(1)=\mathrm{SL}(2,\mathbb{Z})$ is the full modular group. Clearly 
$\Gamma(N) \subset \Gamma_1(N)\subset \Gamma_0(N) \subset \Gamma(1)$, moreover $\Gamma(N)$, the \emph{principal 
congruence subgroup} of level $N$, is the kernel of the homomorphism $\mathrm{SL}(2,\mathbb{Z}) \rightarrow \mathrm{SL}(2,
\mathbb{Z}/N\mathbb{Z})$. Let $\Gamma$ be a discrete subgroup $\Gamma(N) \subset \Gamma \subset \Gamma(1)$ such that 
$N$ is the smallest such integer. A \emph{fundamental domain} is the closure of $\mathbb{H}/\Gamma$, e.g., $\mathbb{H}/\Gamma(1)
= \left\{z\in\mathbb{H} \,|\, -\frac{1}{2} < \mathrm{Re }(z) < \frac{1}{2}, |z| > 1 \right\}$. 

A modular form $f : \mathbb{H}\to\mathbb{C}$ of \emph{weight} $k$ and \emph{level} $N$ is a holomorphic form on the upper 
half plane $\mathbb{H}$ that transforms as
\begin{equation}
f\left(\frac{az + b}{cz + d}\right) = \chi_N(d) (cz + d)^k f(z), \quad \text{where } 
\gamma = \left(\begin{array}{cc}
a & b\\ 
c & d\end{array}\right) \in \Gamma 
\label{ModTrans}
\end{equation}
under the action of a discrete subgroup $\Gamma(N) \subset \Gamma \subset \Gamma(1)$. 

A modular form of SL(2,$\mathbb{Z}$) (i.e., level 1) admits a Fourier expansion ($q$-expansion)  in terms of $q=e^{2\pi i z}$ as
\begin{equation}
f = \sum_{n=0}^\infty a(n) q^n = \sum_{n=0}^\infty a(n) e^{2\pi i n z} \label{qExpand} 
\end{equation}
It is called a \emph{cusp form} if $f$ vanishes as $\mathrm{Im }(z)\to i\infty$, or equivalently at $q=0$. This requires $a(0)=0$. 
It is conventional, and convenient for many purposes, to normalise such that $a(1)=1$. We shall assume this in what follows.
These forms are related to scaling functions on $\mathbb{C}/\Lambda$ where $\Lambda$ is a lattice left invariant by the action
of the modular group. The set of modular forms of weight $k$ form a \emph{finite dimensional} complex vector space 
$M_k\left(\Gamma(1)\right)$ and the set of cusp forms of weight $k$ is a subspace $S_k\left(\Gamma(1)\right)$. Similar notions 
exist for modular forms of subgroups $\Gamma \subset \Gamma(1)$ of level $N$. The cusp forms of a more general modular 
group $\Gamma$ vanish as $z$ approaches certain rational points on $\mathbb{R} = \partial\mathbb{H}$. In the fundamental 
domain these are images of  $\mathrm{Im }(z)\to i\infty$ under $\Gamma(1)/\Gamma(N)$.  

The Dirichlet series associated to a cusp form $f$ is  
\begin{equation}
L(s,f) =  \sum_{n=1}^\infty \frac{a(n)}{n^s} = 1 + \frac{a(2)}{2^{s}} + \frac{a(3)}{3^{s}} + \frac{a(4)}{4^{s}} + \cdots
\label{ModularL}
\end{equation}
where $s\in\mathbb{C}$ and the normalisation is such that $a(1)=1$. The series converges uniformly to a holomorphic function 
of $s$ in the region $\mathrm{Re }(s) > \sigma$ as long as the coefficients $a(n)$ are bounded by some power $n^\sigma$. The 
corresponding $L$-function associated to the cusp form $f$ is then defined by analytic continuation to the complex $s$-plane. If 
$f$ is a cusp form of weight $k$, then the series above converges in $\mathrm{Re }(s) > 1 + \frac{k}{2}$.  The series can also be
expressed as
\begin{equation}
L(s,f) = \mathcal{M}[f](s) = \frac{(2\pi)^s}{\Gamma(s)} \int_0^\infty dy\, y^{s-1} f(iy)  \label{FMellinL}
\end{equation}
i.e., as a Mellin transform of $f(iy)$.
  
Remarkably, the $L$-function of a cusp form $f$ (or modular $L$-function for brevity) \cref{ModularL} has an Euler product form
\begin{equation}
L(s,f) = \sum_{n=1}^\infty \frac{a(n)}{n^s}  = \prod_{p\,\in\,{\mathrm{primes}}} \frac{1}{\left(1 - a(p) p^{-s} +\chi(p) 
p^{k-1} p^{-2s}\right)} \label{EulerModL}   
\end{equation}  
The coefficients $a(n)$ have very interesting properties, to which we shall return in a moment. As an example, consider one of the 
well known modular $L$-functions related to the discriminant function 
$\Delta(z) = \left(\sqrt{2\pi}\,\eta(z)\right)^{24} = (2\pi)^{12}\, q\, \displaystyle{\prod_{n=1}^\infty \left(1-q^n\right)}$, 
where $\eta(z)$ is the Dedekind $\eta$-function. It is a holomorphic modular form of weight 12 (and level 1) that vanishes as 
$z\to i\infty$, i.e., it is a cusp form. The coefficients in its $q$-expansion 
\begin{equation}
\Delta(z)  =  \sum_{n=1}^\infty \tau(n) q^n
\end{equation}
define the function $\tau : \mathbb{N}\rightarrow \mathbb{Z}$, known as the Ramanujan $\tau$-function. They satisfy the following
properties
\begin{align}
\begin{split}
\tau(mn) &= \tau(m) \tau(n) \; \text{if gcd }(m,n)=1\\
\tau(p^{m+1}) &= \tau(p) \tau(p^m) - p^{11} \tau(p^{m-1}) \;\, \text{for } p \:\; \text{a prime and } m>0 \\
|\tau(p)| &\le 2 p^{11/2}
\end{split}\label{RamanujanTau}
\end{align}
as were conjectured by Ramanujan. The first two were proved by Mordell soon after the conjecture, while the bound was proved by 
Deligne after many years. More generally, the coefficients $a(n)$ of a modular form of weight $k$ and level $N$ satisfy
\begin{align}
\begin{split}
a(mn) &= a(m) a(n) \; \text{if gcd }(m,n)=1\\
a(p^{m+1}) &= a(p)\, a(p^m) - \chi(p) p^{k-1} a(p^{m-1}) \;\, \text{for } p \:\; \text{a prime and } m>0 
\end{split}\label{Recursion}
\end{align}
as follows from comparing the series and product formulas in \cref{EulerModL}. Therefore, the coefficients $a(n)$ define a 
multiplicative character \cite{Serre:Course} on $\mathbb{N}$. The convergence of the series also puts a bound on the growth of 
the coefficients.

\bigskip

Let us consider a prime factor in \cref{EulerModL} for a fixed $p$, and call it a \emph{local} modular $L$-function at a prime $p$ in
analogy with the local zeta function \cref{locZTrace} at $p$, and factorise it as follows  
\begin{equation}
L_{(p)}(s,f) =  \frac{1}{\left(1 - a(p) p^{-s} +\chi(p) p^{k-1} p^{-2s}\right)} = \frac{1}{\left(1 - a_1(p) p^{-s}\right) 
\left(1 - a_2(p) p^{-s}\right)} \label{LocModL}
\end{equation}
Evidently, consistency demands that the sum and the product of the coefficient functions of the factorised form gives
\begin{equation}
a_1(p) + a_2(p) = a(p) \quad\text{and}\quad a_1(p)\, a_2(p) = \chi(p) p^{k-1} \label{Consistency}
\end{equation}
Notice that $a_1(p)$ and $a_2(p)$ are only defined for the (fixed) prime $p$. For the positive powers of $p$ we take 
$a_i : p^\mathbb{N} \rightarrow \mathbb{C}$ by $a_i(p^n) = \left(a_i(p)\right)^n$. 

In the region of convergence, we can realise \cref{LocModL} as infinite geometric series
\begin{align}
L_{(p)}(s,f) &= \sum_{m=0}^\infty \sum_{\ell=0}^m (-1)^\ell {m\choose\ell} \left(a(p)\right)^{m-\ell} \left(\chi(p) p^{k-1}\right)^\ell
p^{-(m+\ell)s} \nonumber\\
&= \sum_{m_1=0}^\infty \sum_{m_2=0}^\infty \left(a_1(p)\right)^{m_1} \left(a_2(p)\right)^{m_2} p^{-(m_1+m_2)s}  
\label{ModLGSeries}
\end{align}
Hence, by comparing terms of powers of $p^{-s}$, we find
\begin{equation}
\sum_{n=0}^m \left(a_1(p)\right)^{m-n} \left(a_2(p)\right)^n = \sum_{\ell=0}^{\left\lfloor{\frac{m}{2}}\right\rfloor} (-1)^\ell 
{m-\ell\choose\ell} \left(a(p)\right)^{m-2\ell} \left(\chi(p) p^{k-1}\right)^\ell \label{RelationAA1A2}
\end{equation}
which can be checked to be consistent with \cref{Consistency}.

\bigskip

{}From here onwards, we can closely follow the steps in the extension of the Dirichlet character in \cref{sec:DChar}, and define the 
functions $\mathfrak{a}_1$ and $\mathfrak{a}_2$ as follows: 
\begin{equation}
\text{For } i=1,2,\;\mathfrak{a}_i(p^{n}) = \left\{\begin{array}{ll}
 \left(a_i(p)\right)^{n} & \text{ if } \: a_i(p) \ne 0\:\text { and }\, n\ge 0\\
 \left(a_i(p)\right)^{n} = \displaystyle{\frac{1}{\left(a_i(p)\right)^{|n|}}} & \text{ if } \: a_i(p) \ne 0\:\text { and }\, n < 0\\
 0 & \text{ if } \: a_i(p) = 0 \end{array} 
  \right. \label{AChar}
\end{equation}
The extended functions are $\mathfrak{a}_i : p^{\mathbb{Z}} \rightarrow \mathbb{C}$. With the help of these functions, we define
the gamma functions $\Gamma_{\mathfrak{a}_1}(s,f)$ and $\Gamma_{\mathfrak{a}_2}(s,f)$, related to the modular function $f$ as
\begin{align}
\Gamma_{(\mathfrak{a}_i,f)} (s) 
&= \int_{\Qp} d\xi \, e^{2\pi i\xi} |\xi|_p^{s-1}\, \mathfrak{a}_i\!\left(|\xi|_p^{-1}\right) , \quad \text{for } i=1,2\nonumber\\
&= \frac{a_i(p) -  p^{s-1}}{a_i(p)\left(1 - a_i(p) p^{-s}\right) }   \label{GammA}
\end{align}
for non-vanishing $a_i(p)$, otherwise for $a_i(p)=0$, we take the gamma function to be zero. Finally, let us define a pair of 
generalised Vladimirov derivatives associated with the modular function $f$, as  
\begin{equation}\label{AVladD}
D_{(\mathfrak{a}_i,f)}^\alpha g(\xi) = \frac{1}{\Gamma_{(\mathfrak{a}_i,f)}(-\alpha)}\, \int_{\Qp} d\xi'\, 
\frac{g(\xi') - g(\xi)}{|\xi'-\xi|_p^{\alpha+1}}\, \mathfrak{a}_i(|\xi'-\xi|_p^{-1}) 
\end{equation}
unless, $\mathfrak{a}_i=0$, in which case we define $D_{(\mathfrak{a}_i,f)} = 1$. Repeating our previous arguments {\em
mutatis mutandis}, it is straightforward to check that the Kozyrev wavelets are eigenfunctions of the derivative operators
\cref{AVladD}
\begin{equation}
D_{(\mathfrak{a}_i,f)}^\alpha \psi_{n,m,j}^{(p)} (\xi) = \left(a_i(p) p^{\alpha}\right)^{(1-n)} \psi_{n,m,j}^{(p)} (\xi),\quad
\text{for }i=1,2 \label{AVladOnKozy}
\end{equation}
with eigenvalues involving the factorised coefficients. In proving this, one needs to sum over an infinite geometric series which
converges for $\mathrm{Re }(s) > 1 + \frac{k}{2}$.
 
The equations above, together with \cref{LocModL} and \cref{ModLGSeries}, allow us to write the local factors of the modular 
$L$-function associated $f$ as
\begin{eqnarray}
L_{(p)}(s,f) &=& \left(\sum_{m_1=0}^\infty \left(a_1(p) p^{-s}\right)^{m_1}\right) \left(\sum_{m_2=0}^\infty 
\left(a_2(p) p^{-s}\right)^{m_2}\right) \nonumber\\ 
&=& \sum_{m=0}^\infty \left(\sum_{m_1=0}^m \left(a_1(p)\right)^{m_1} \left(a_2(p)\right)^{m-m_1}\right) p^{-sm}\nonumber\\ 
&=& \mathrm{Tr}_{\mathcal{H}_-\otimes\mathcal{H}_-} \left(D_{(\mathfrak{a}_1,f)}^{-s} \otimes D_{(\mathfrak{a}_2,f)}^{-s}\right)  
\label{LocModLTrace}
\end{eqnarray}
The trace is now over a subspace $\mathcal{H}_-\otimes\mathcal{H}_-$ of the tensor product space $L^2(\Qp) \otimes L^2(\Qp)$. The 
sum in the trace is over the wavelets with negative scaling quantum number, hence compactly supported in $\Zp^2 \subset \Qp^2$. 
Therefore, the double sum  above span the upper right quadrangle of the lattice $\mathbb{Z}^2$ (see \cref{fig:lattice}).

\begin{figure}[h]
\centering
\begin{tikzpicture}[scale=0.7,every node/.style={scale=0.8}]
\filldraw[line width=1, blue] (0,3.5) -- (3.5,0);
\fill[lightgray,opacity=0.3]    (3.5,-0.15) -- (5.25,-0.15) -- (5.25,5.25) -- (-0.15,5.25) -- (-0.15,3.5) -- cycle;
\node[align=center,mark size=2pt,color=black] at (0,0) {\pgfuseplotmark{*}};
\node[mark size=2pt,color=black] at (0.5,0) {\pgfuseplotmark{*}};
\node[mark size=2pt,color=black] at (1,0) {\pgfuseplotmark{*}};
\node[mark size=2pt,color=black] at (1.5,0) {\pgfuseplotmark{*}};
\node[mark size=2pt,color=black] at (2,0) {\pgfuseplotmark{*}};
\node[mark size=2pt,color=black] at (2.5,0) {\pgfuseplotmark{*}};
\node[mark size=2pt,color=black] at (3,0) {\pgfuseplotmark{*}};
\node[mark size=2pt,color=black] at (3.5,0) {\pgfuseplotmark{*}};
\node[mark size=2pt,color=black] at (4,0) {\pgfuseplotmark{*}};
\node[mark size=2pt,color=black] at (4.5,0) {\pgfuseplotmark{*}};
\node[mark size=2pt,color=black] at (5,0) {\pgfuseplotmark{*}};
\node[align=center,mark size=2pt,color=black] at (0,0.5) {\pgfuseplotmark{*}};
\node[mark size=2pt,color=black] at (0.5,0.5) {\pgfuseplotmark{*}};
\node[mark size=2pt,color=black] at (1.0,0.5) {\pgfuseplotmark{*}};
\node[mark size=2pt,color=black] at (1.5,0.5) {\pgfuseplotmark{*}};
\node[mark size=2pt,color=black] at (2,0.5) {\pgfuseplotmark{*}};
\node[mark size=2pt,color=black] at (2.5,0.5) {\pgfuseplotmark{*}};
\node[mark size=2pt,color=black] at (3,0.5) {\pgfuseplotmark{*}};
\node[mark size=2pt,color=black] at (3.5,0.5) {\pgfuseplotmark{*}};
\node[mark size=2pt,color=black] at (4,0.5) {\pgfuseplotmark{*}};
\node[mark size=2pt,color=black] at (4.5,0.5) {\pgfuseplotmark{*}};
\node[mark size=2pt,color=black] at (5,0.5) {\pgfuseplotmark{*}};
\node[mark size=2pt,color=black] at (0,1) {\pgfuseplotmark{*}};
\node[mark size=2pt,color=black] at (0.5,1) {\pgfuseplotmark{*}};
\node[mark size=2pt,color=black] at (1,1) {\pgfuseplotmark{*}};
\node[mark size=2pt,color=black] at (1.5,1) {\pgfuseplotmark{*}};
\node[mark size=2pt,color=black] at (2,1) {\pgfuseplotmark{*}};
\node[mark size=2pt,color=black] at (2.5,1) {\pgfuseplotmark{*}};
\node[mark size=2pt,color=black] at (3,1) {\pgfuseplotmark{*}};
\node[mark size=2pt,color=black] at (3.5,1) {\pgfuseplotmark{*}};
\node[mark size=2pt,color=black] at (4,1) {\pgfuseplotmark{*}};
\node[mark size=2pt,color=black] at (4.5,1) {\pgfuseplotmark{*}};
\node[mark size=2pt,color=black] at (5,1) {\pgfuseplotmark{*}};
\node[mark size=2pt,color=black] at (0,1.5) {\pgfuseplotmark{*}};
\node[mark size=2pt,color=black] at (0.5,1.5) {\pgfuseplotmark{*}};
\node[mark size=2pt,color=black] at (1.0,1.5) {\pgfuseplotmark{*}};
\node[mark size=2pt,color=black] at (1.5,1.5) {\pgfuseplotmark{*}};
\node[mark size=2pt,color=black] at (2,1.5) {\pgfuseplotmark{*}};
\node[mark size=2pt,color=black] at (2.5,1.5) {\pgfuseplotmark{*}};
\node[mark size=2pt,color=black] at (3,1.5) {\pgfuseplotmark{*}};
\node[mark size=2pt,color=black] at (3.5,1.5) {\pgfuseplotmark{*}};
\node[mark size=2pt,color=black] at (4,1.5) {\pgfuseplotmark{*}};
\node[mark size=2pt,color=black] at (4.5,1.5) {\pgfuseplotmark{*}};
\node[mark size=2pt,color=black] at (5,1.5) {\pgfuseplotmark{*}};
\node[mark size=2pt,color=black] at (0,2) {\pgfuseplotmark{*}};
\node[mark size=2pt,color=black] at (0.5,2) {\pgfuseplotmark{*}};
\node[mark size=2pt,color=black] at (1,2) {\pgfuseplotmark{*}};
\node[mark size=2pt,color=black] at (1.5,2) {\pgfuseplotmark{*}};
\node[mark size=2pt,color=black] at (2,2) {\pgfuseplotmark{*}};
\node[mark size=2pt,color=black] at (2.5,2) {\pgfuseplotmark{*}};
\node[mark size=2pt,color=black] at (3,2) {\pgfuseplotmark{*}};
\node[mark size=2pt,color=black] at (3.5,2) {\pgfuseplotmark{*}};
\node[mark size=2pt,color=black] at (4,2) {\pgfuseplotmark{*}};
\node[mark size=2pt,color=black] at (4.5,2) {\pgfuseplotmark{*}};
\node[mark size=2pt,color=black] at (5,2) {\pgfuseplotmark{*}};
\node[mark size=2pt,color=black] at (0,2.5) {\pgfuseplotmark{*}};
\node[mark size=2pt,color=black] at (0.5,2.5) {\pgfuseplotmark{*}};
\node[mark size=2pt,color=black] at (1.0,2.5) {\pgfuseplotmark{*}};
\node[mark size=2pt,color=black] at (1.5,2.5) {\pgfuseplotmark{*}};
\node[mark size=2pt,color=black] at (2,2.5) {\pgfuseplotmark{*}};
\node[mark size=2pt,color=black] at (2.5,2.5) {\pgfuseplotmark{*}};
\node[mark size=2pt,color=black] at (3,2.5) {\pgfuseplotmark{*}};
\node[mark size=2pt,color=black] at (3.5,2.5) {\pgfuseplotmark{*}};
\node[mark size=2pt,color=black] at (4,2.5) {\pgfuseplotmark{*}};
\node[mark size=2pt,color=black] at (4.5,2.5) {\pgfuseplotmark{*}};
\node[mark size=2pt,color=black] at (5,2.5) {\pgfuseplotmark{*}};
\node[mark size=2pt,color=black] at (0,3) {\pgfuseplotmark{*}};
\node[mark size=2pt,color=black] at (0.5,3) {\pgfuseplotmark{*}};
\node[mark size=2pt,color=black] at (1,3) {\pgfuseplotmark{*}};
\node[mark size=2pt,color=black] at (1.5,3) {\pgfuseplotmark{*}};
\node[mark size=2pt,color=black] at (2,3) {\pgfuseplotmark{*}};
\node[mark size=2pt,color=black] at (2.5,3) {\pgfuseplotmark{*}};
\node[mark size=2pt,color=black] at (3,3) {\pgfuseplotmark{*}};
\node[mark size=2pt,color=black] at (3.5,3) {\pgfuseplotmark{*}};
\node[mark size=2pt,color=black] at (4,3) {\pgfuseplotmark{*}};
\node[mark size=2pt,color=black] at (4.5,3) {\pgfuseplotmark{*}};
\node[mark size=2pt,color=black] at (5,3) {\pgfuseplotmark{*}};
\node[mark size=2pt,color=black] at (0,3.5) {\pgfuseplotmark{*}};
\node[mark size=2pt,color=black] at (0.5,3.5) {\pgfuseplotmark{*}};
\node[mark size=2pt,color=black] at (1.0,3.5) {\pgfuseplotmark{*}};
\node[mark size=2pt,color=black] at (1.5,3.5) {\pgfuseplotmark{*}};
\node[mark size=2pt,color=black] at (2,3.5) {\pgfuseplotmark{*}};
\node[mark size=2pt,color=black] at (2.5,3.5) {\pgfuseplotmark{*}};
\node[mark size=2pt,color=black] at (3,3.5) {\pgfuseplotmark{*}};
\node[mark size=2pt,color=black] at (3.5,3.5) {\pgfuseplotmark{*}};
\node[mark size=2pt,color=black] at (4,3.5) {\pgfuseplotmark{*}};
\node[mark size=2pt,color=black] at (4.5,3.5) {\pgfuseplotmark{*}};
\node[mark size=2pt,color=black] at (5,3.5) {\pgfuseplotmark{*}};
\node[mark size=2pt,color=black] at (0,4) {\pgfuseplotmark{*}};
\node[mark size=2pt,color=black] at (0.5,4) {\pgfuseplotmark{*}};
\node[mark size=2pt,color=black] at (1,4) {\pgfuseplotmark{*}};
\node[mark size=2pt,color=black] at (1.5,4) {\pgfuseplotmark{*}};
\node[mark size=2pt,color=black] at (2,4) {\pgfuseplotmark{*}};
\node[mark size=2pt,color=black] at (2.5,4) {\pgfuseplotmark{*}};
\node[mark size=2pt,color=black] at (3,4) {\pgfuseplotmark{*}};
\node[mark size=2pt,color=black] at (3.5,4) {\pgfuseplotmark{*}};
\node[mark size=2pt,color=black] at (4,4) {\pgfuseplotmark{*}};
\node[mark size=2pt,color=black] at (4.5,4) {\pgfuseplotmark{*}};
\node[mark size=2pt,color=black] at (5,4) {\pgfuseplotmark{*}};
\node[mark size=2pt,color=black] at (0,4.5) {\pgfuseplotmark{*}};
\node[mark size=2pt,color=black] at (0.5,4.5) {\pgfuseplotmark{*}};
\node[mark size=2pt,color=black] at (1.0,4.5) {\pgfuseplotmark{*}};
\node[mark size=2pt,color=black] at (1.5,4.5) {\pgfuseplotmark{*}};
\node[mark size=2pt,color=black] at (2,4.5) {\pgfuseplotmark{*}};
\node[mark size=2pt,color=black] at (2.5,4.5) {\pgfuseplotmark{*}};
\node[mark size=2pt,color=black] at (3,4.5) {\pgfuseplotmark{*}};
\node[mark size=2pt,color=black] at (3.5,4.5) {\pgfuseplotmark{*}};
\node[mark size=2pt,color=black] at (4,4.5) {\pgfuseplotmark{*}};
\node[mark size=2pt,color=black] at (4.5,4.5) {\pgfuseplotmark{*}};
\node[mark size=2pt,color=black] at (5,4.5) {\pgfuseplotmark{*}};
\node[mark size=2pt,color=black] at (0,5) {\pgfuseplotmark{*}};
\node[mark size=2pt,color=black] at (0.5,5) {\pgfuseplotmark{*}};
\node[mark size=2pt,color=black] at (1,5) {\pgfuseplotmark{*}};
\node[mark size=2pt,color=black] at (1.5,5) {\pgfuseplotmark{*}};
\node[mark size=2pt,color=black] at (2,5) {\pgfuseplotmark{*}};
\node[mark size=2pt,color=black] at (2.5,5) {\pgfuseplotmark{*}};
\node[mark size=2pt,color=black] at (3,5) {\pgfuseplotmark{*}};
\node[mark size=2pt,color=black] at (3.5,5) {\pgfuseplotmark{*}};
\node[mark size=2pt,color=black] at (4,5) {\pgfuseplotmark{*}};
\node[mark size=2pt,color=black] at (4.5,5) {\pgfuseplotmark{*}};
\node[mark size=2pt,color=black] at (5,5) {\pgfuseplotmark{*}};
\end{tikzpicture}
\caption{{\small The sum in Eq.(\ref{LocModLTrace}) is over the upper quadrangle of the lattice $\mathbb{Z}^2$, 
while the sum in Eq.(\ref{HeckeLike}) is over the region above and right of the blue line. }}
\label{fig:lattice}
\end{figure}

It is interesting to note that if we conjugate the operators \cref{AVladD} by the following combination of the raising/lowering 
operators in \cref{raise:lower} and then take the trace\footnote{Similarly, although not of much interest, the operators in 
\cref{locZTrace,locLTrace} may be conjugated to get $\mathrm{Tr}_{\mathcal{H}_-} \left(a_+^\ell D^{-s}_{\mathfrak{x}} 
a_-^\ell\right) = \left(\chi(p) p^{-s}\right)^\ell L_p(s,\chi)$, i.e, the $\ell$-th coefficient as a prefactor.}, we get
\begin{eqnarray}
\mathrm{Tr} \left(\sum_{k=0}^\ell a_+^k D_{(\mathfrak{a}_1,f)}^{-s} a_-^k \otimes a_+^{\ell-k} 
D_{(\mathfrak{a}_2,f)}^{-s} a_-^{\ell-k}\right) 
&=& \displaystyle{\sum_{m=\ell}^\infty \left(\sum_{m_1=0}^m \left(a_1(p)\right)^{m_1} \left(a_2(p)\right)^{m-m_1}\right) 
p^{-sm}}\nonumber\\
&=& \displaystyle{\left(\sum_{k=0}^\ell \left(a_1(p)\right)^{k} \left(a_2(p)\right)^{\ell - k}\right) p^{-s\ell}} L_{(p)}(s,f) \nonumber\\
&=& a(p^\ell) p^{-s\ell} L_{(p)}(s,f) \label{HeckeLike}
\end{eqnarray}
In the above, the raising and lowering operators refer to their restrictions to the subspace $\mathcal{H}_-$, however, our use 
of the same notation will hopefully not cause any confusion. Therefore, the action of $a_+$, or more accurately that of
$a_+\big|_{\mathcal{H}_-}$, on the `ground state' $|\Omega\rangle$ corresponding to $n=1$ (mother wavelet with scaling
quantum number zero) is to annihilate it, i.e., $a_+ |\Omega\rangle = 0$. The effect of the sum of the conjugated operators is 
to sum over a part of the lattice $\mathbb{Z}^2$ satisfying $m_1+m+2\ge \ell$ (see \cref{fig:lattice}).  This action mimics the 
action of the Hecke operator $T_{p^\ell}$ on the $L$-function. 

Recall that the Hecke operators $T(m)$, $m\in\mathbb{N}$, are a set of commuting operators whose action on the modular
form is to return the coefficients in the $q$-expansion as eigenvalues
\begin{equation}
T(m) f(z) = a(m) f(z) \label{HeckeOnModForm}
\end{equation}
In other words a modular form is an eigenvector of the Hecke operators with the eigenvalues as the coefficients in its 
$q$-expansion. They satisfy
\begin{align}
\begin{split}
T(m) T(n) &= T(mn) \text{ for } m\nmid n\\
T(p) T(p^\ell) &= T(p^{\ell+1}) + \chi(p) p^{k-1} T(p^{\ell-1}) 
\end{split}
\label{HeckeAlg}
\end{align}
A modular forms can also be understood\cite{Serre:Course,Koblitz:ECMF} as a sum over the nodes of a lattice $\Lambda$ in 
$\mathbb{C}$. {}From this point of view, the action of the Hecke operator $T(n)$ involves sublattices $\Lambda' \subset \Lambda$ 
of index $n$: $T(n) \Lambda =\! \displaystyle{\sum_{[\Lambda:\Lambda']=n}}\! \Lambda'$.   

\section{Summary}
We have generalised the notion of the Vladimirov derivative, a pseudodifferential operator on Bruhat-Schwartz functions on 
the $p$-adic space $\Qp$, by including the \emph{twist} by multiplicative characters. The simplest of the multiplicative characters
in this context are Dirichlet characters. We show that the wavelet functions, defined by Kozyrev, are eigenfunctions of the twisted
operators with eigenvalues that depend on the character and the scaling properties of the wavelets. This allowed us to write the
Dirichlet $L$-series (including the Riemann zeta function, which corresponds to the trivial character) as traces of appropriate
twisted Vladimirov derivatives on a subspace of Bruhat-Schwartz functions supported on the compact subset of $p$-adic integers.
Along the way, we have defined a generalised class of Gelfand-Graev-Tate gamma functions on $\Qp$ corresponding to twist by 
multiplicative characters. In \cref{sec:VladDLFn}, we have further generalised our construction to the $L$-series associated to the 
cusp forms of congruence subgroups of the modular group SL(2,$\mathbb{Z}$). It would be interesting to realise these $L$-functions
as the partition function of a `statistical system'. We hope to report on the last point in the near future.

\bigskip

\noindent{\bf Acknowledgments:} One of us (DG) presented the results in this paper in the \emph{VII-th International Conference 
$p$-Adic Mathematical Physics \&\ its Applications} held at the Universidade Beira Interior, Covilh\~{a}, Portugal during Sep 30--Oct 
4, 2019, as well as in the \emph{National String Meeting 2019} held at IISER Bhopal, India during Dec 22--27, 2019. We thank the 
participants of these meetings, especially P.~Bradley and W.~Zu\~{n}iga-Galindo for their comments. We also thank V.~Patankar for 
discussions. 
 
\bigskip

\noindent{\bf Note added:} Recently a paper (\href{https://arxiv.org/abs/2001.01721}{arXiv:2001.01721 [hep-th]}) that discussed the 
construction of pseudodifferential operators of similar type appeared on the arXiv.

\newpage

\bibliographystyle{hieeetr}
\bibliography{PseuD}{}

\end{document}